\documentclass[amstex,16pt, amssymb]{article}

\usepackage{mathtext}
\usepackage[cp1251]{inputenc}
\usepackage[T2A]{fontenc}
\usepackage[dvips]{graphicx}
\usepackage{amsmath}
\usepackage{amssymb}
\usepackage{amsxtra}
\usepackage{latexsym}
\usepackage{ifthen}
\usepackage{amsthm}

\numberwithin{equation}{section}

\begin{document}

\title{Toward the theory of semi-linear
equations}

\author{Vladimir Gutlyanski\u\i{}, Olga Nesmelova, Vladimir Ryazanov}

\date{}

\maketitle 

\begin{abstract}
In this paper we study the semilinear partial differential
equations in the plane, the linear part of which is written in a
divergence form. The main result is given as a factorization
theorem. This theorem states that every  weak solution of such an
equation can be represented as a composition of a weak solution of
the corresponding isotropic equation in a canonical domain and a
quasiconformal mapping agreed with a matrix-valued measurable
coefficient appearing in the divergence part of the equation. The
latter makes it possible, in particular, to remove the regularity
restrictions on the boundary in the study of boundary value problems
for such semilinear equations.
\end{abstract}



\normalsize \baselineskip=18.5pt

\vskip 1cm


\section{Introduction}

It is well known that the first order linear partial differential
equation \begin{equation}
\label{Beltrami} \omega_{\bar z} =\mu(z)\,
\omega_z,\,\,\,\,z\in\Omega\subseteq{\Bbb C}, 
\end{equation} where
$\omega_z={1\over 2}(\omega_x-i\omega_y),$ $\omega_{\bar z}={1\over
2}(\omega_x+i\omega_y),$ $z=x+i\, y,$ which has come to be known as
the Beltrami equation, turned out to be instrumental in the study of
Riemann surfaces, Teichm\"uller spaces, Kleinian groups, meromorphic
functions, low dimensional topology, holomorphic motion complex
dynamics, Clifford analysis and control theory. As known, a
$K-$quasiconformal mapping
 $\omega:\Omega\to{\Bbb C},$ $K\geq 1,$ is just an orientation-preserving
homeomorphic $W_{\rm loc}^{1,2}(\Omega)$ solution to  the Beltrami
equation when the measurable coefficient $\mu$ satisfies the strong
ellipticity condition $|\mu(z)|\leq (K-1)/(K+1)$ almost everywhere
in $\Omega.$ Note in particular, that if $\mu=0$ in a domain
$\Omega\subset{\Bbb C},$ then the Beltrami equation reduces to the
Cauchy--Riemann equation and a solution $\omega$ is analytic in
$\Omega,$ see e.g. \cite{Ahlfors:book}, \cite{LV:book}, see also
\cite{GRSY}, and the references therein.
\par
We also would like to pay attention to a strong interaction between
linear and non-linear elliptic systems in the plane and
quasiconformal mappings. The most general first order linear
homogeneous elliptic system with real coefficients can be written in
the form $\omega_{\bar z}+\mu(z)\, \omega_z +\nu(z)\,
\overline{\omega_z}=0,$ with measurable coefficients $\mu$ and $\nu$
such that $|\mu|+|\nu|\leq (K-1)/(K+1)<1.$ This equation is a
particular case of a non-linear first order system $\omega_{\bar
z}=H(z,\omega_z)$ where $H:G\times{\Bbb C}\to{\Bbb C}$ is Lipschitz
in the second variable,
$$|H(z,w_1)-H(z,w_2)|\ \leq\  {K-1\over K+1}\ |w_1-w_2|\ ,\,\,\,\, H(z,0)\equiv 0\ .$$
The principal feature of the above equation is that the difference
of two solutions  need not solve the same equation but each
solution can be represented as a composition of a quasiconformal
homeomorphism and an analytic function. Thus quasiconformal
mappings become the central tool for the study of these non-linear
systems. A rather comprehensive treatment of the present state of
the theory is given in the excellent book of Astala, Iwaniec and
Martin \cite{AIM-book}. This book contains also an exhaustive
bibliography on the topic.
 In particular, the following fundamental Harmonic
Factorization Theorem for the uniformly elliptic divergence
equations \begin{equation}
\label{de} {\rm div\,}A(z,\nabla
u)=0,\,\,\,z\in\Omega, 
\end{equation} holds, see \cite{AIM-book}, Theorem
16.2.1: Every solution $u\in W^{1,2}_{\rm loc}(\Omega)$
 of the  equation (\ref{de})
can be represented as $u(z) = h(\omega(z)),$ where $\omega :
\Omega\to G$ is $K-$quasiconformal and $h$ is harmonic on $G.$
\par
The main goal of this paper is  to point out another application
of quasiconformal mappings to the study of some {\it nonlinear}
partial differential equations in the plane. Namely, we will deal
with co-called {\it semi-linear} partial differential equations,
linear part of which contains the elliptic operator in the
divergence form ${\rm div}\,[A(z)\nabla u(z)],$ where the
 matrix function $A(z)$ is assumed to be symmetric, measurable
and satisfy the uniform ellipticity condition \begin{equation}
\label{ellipticity}
{1\over K}\, |\xi|^2\ \leq \ \langle A(z)\, \xi,\xi\rangle\ \leq\
K\, |\xi|^2 \,\,\,\mbox{a.e. in $\Omega$} 
\end{equation}
 for every $\xi\in{\Bbb
C},$ where $1\leq K<\infty$. Some examples of such semi-linear
equations include in anisotropic case the nonlinear heat equation
like \begin{equation}\label{he} u_t-{\rm div}\,[A(z)\nabla u(z)]=f(u), 
\end{equation}
 (the
same equation describes the brownian motion, diffusion models of
population dynamics, and many other phenomena), nonlinear wave
equation \begin{equation}
\label{we} u_{tt}-{\rm div}\,[A(z)\nabla u(z)]=f(u), 
\end{equation}
with continuous functions $f$ or nonlinear Schr\"odinger equation
\begin{equation}\label{se} iu_t+{\rm div}\,[A(z)\nabla u(z)]=k|u|^pu, 
\end{equation} as well
as their stationary counterparts. We show that every weak solution
of semi-linear equations of such type can be factorized  as the
composition of a solution to the corresponding isotropic equation
and a quasiconformal mapping agreed with the matrix function $A(z).$
If we deal, for example, with the semi-linear equation \begin{equation}
\label{b1}
{\rm div}\,[A(z)\nabla u(z)]=f(u),\,\,\,\,f\geq 0, 
\end{equation} then we prove
the following generalization of the mentioned above Factorization
theorem: Every solution $u\in W^{1,2}_{\rm loc}(\Omega)$
 of the semi-linear equation (\ref{b1})
can be represented as $u(z) = h(\omega(z)),$ where $\omega :
\Omega\to G$ is a quasiconformal mapping agreed with the matrix
function $A$ and $h$ is subharmonic on $G.$

\par
The paper is organized as follows. In the next section we recall
basic results from the theory of quasiconformal mappings and their
connection with linear elliptic differential equations in the
divergence form. In Section 3  we prove the basic identity
concerning the linear elliptic operator in the divergence form. The
main result of this section is Proposition 3.1. In Section 4 we
prove the principal Factorization Theorem 4.1 and their corollaries
for some elliptic semi-linear model equations in the plane. Some
applications of Factorization Theorem to existence theorems for
elliptic semi-linear equations are given in Section 5. Its
applications to the free boundary problems ("dead zones") can be
found in Section 6. A similar discussion of equations of the heat
type is in Section 7. Finally, Section 8 contains a discussion of
boundary value problems.

\section{QK-maps and divergent equations}

Let $\Omega$  be a domain in the complex plane ${\Bbb C}.$ Denote by
$M^{2\times 2}(\Omega)$ the class of all $2\times 2$ symmetric
matrix function $A(z)=\{a_{jk}(z)\}$ with measurable entries and
${\rm det}\,A(z)=1,$ satisfying the uniform ellipticity condition
\begin{equation}\label{ellipticity} {1\over K}\, |\xi|^2\ \leq\ \langle A(z)\,
\xi,\,\xi\rangle\ \leq\ K\, |\xi|^2 \,\,\, \,\,\, \,\,\, \mbox{a.e.
in}\,\, \Omega 
\end{equation} for every $\xi\in{\Bbb C}$ where $1\leq
K<\infty$.
\par
Given $A\in M^{2\times 2}(\Omega),$ let us first consider the
following second order elliptic homogeneous equation in the
divergent form \begin{equation}
\label{second} {\rm div}\,(A(z)\nabla\,u)\ =0\
\,\,\,\,\, \,\,\,\mbox{a.e. in $\Omega\ .$} 
\end{equation}
 Equation
(\ref{second}) is interpreted in the distributional sense. That is,
a function $u$ is a weak solution to the equation if it has locally
integrable gradient $\nabla u$ with 
\begin{equation}\label{weaks}
\int_\Omega\langle A(z)\nabla u,\,\nabla \varphi\rangle\ =\
0\,\,\,\,\, \,\,\, \forall\varphi\in C^\infty_0(\Omega)\ . 
\end{equation} This
is meaningful at least for $u\in W^{1,1}_{\rm loc}(\Omega),$ where
$W^{1,p}_{\rm loc}(\Omega)$ stands for the well-known Sobolev space.
Here we will assume a little more regularity, namely that $u\in
C\cap W^{1,2}_{\rm loc}(\Omega),$ because in this paper we will be
dealing with quasiconformal mappings generated by $A$ with the
uniform ellipticity condition (\ref{ellipticity}). We plan to study
also semi-linear {\it degenerate} partial differential equations in
subsequent works.  In that case we will consider the weaker
solutions $u\in C\cap W^{1,1}_{\rm loc}(\Omega).$

\par
Let $A\in M^{2\times 2}(\Omega)$ and $u\in C\cap W^{1,2}_{\rm
loc}(\Omega)$ be a weak solution to (\ref{second}). Then there
exists $v\in C\cap W^{1,2}_{\rm loc}(\Omega),$ called the {\it
stream function} of $u,$ such that 
\begin{equation}
\nabla\,v\ =\ H\,A\, \nabla\,u \,\,\,\,\,\mbox{a.e in $\Omega,$}\,\,\,\mbox{where}\ H\ =\ \left(\begin{array}{ccc} 0 & -1 \\
                            1          & 0 \end{array}\right)\ .
\end{equation} Note by the way that
$$
H^2\ =\ -\ I\ ,\,\,\,\mbox{where}\  I\ =\ \left(\begin{array}{ccc} 1 & 0 \\
                            0         & 1 \end{array}\right)\ ,
$$
i.e. $H$ takes a part of the imaginary unit in the space of $2\times
2$ matrices.

Setting $\omega(z)=u(z)+i\, v(z),$ one writes in complex notations
that $\omega$ satisfies the Beltrami equation 
\begin{equation}\label{Beltrami2}
\omega_{\bar z}(z)\ =\ \mu(z)\ \omega_z(z)\,\,\,\,\,\mbox{a.e. in
$\Omega\ ,$}\end{equation}
where the complex dilatation $\mu(z)$ is given by
\begin{equation}
\label{mu} \mu(z)\ =\ \frac{a_{22}(z)-a_{11}(z)-2i
a_{12}(z)}{{\rm det}\,\,(I+A(z))}\ . 
\end{equation}
 The condition of
ellipticity (\ref{ellipticity}) now is written as
\begin{equation}\label{ellipticity1}
|\mu(z)|\ \leq\ {K-1\over K+1} \,\,\,\,\,\mbox{a.e. in $\Omega\ .$}
\end{equation} And vice versa, given measurable complex valued function $\mu,$
satisfying (\ref{ellipticity1}), one can invert the algebraic system
(\ref{mu}) to obtain that, see e.g. Theorem 16.1.6 in
\cite{AIM-book}, 
\begin{equation}\label{matrix}
A(z)\ =\ \left(\begin{array}{ccc} {|1-\mu|^2\over 1-|\mu|^2}  & {-2{\rm Im}\,\mu\over 1-|\mu|^2} \\
                            {-2{\rm Im}\,\mu\over 1-|\mu|^2}          & {|1+\mu|^2\over 1-|\mu|^2}
                              \end{array}\right).
\end{equation}

Thus, given any $A\in M^{2\times 2}(\Omega),$ one produces by
(\ref{mu}) the complex dilatation $\mu(z)$ for which in turn, by the
Measurable  Riemann mapping theorem,  see e.g. Theorem V.B.3 in
\cite{Ahlfors:book} and Theorem V.1.3 in \cite{LV:book},
 the Beltrami equation (\ref{Beltrami2})
generates as its solution a quasiconfomal homeomorphism
$\omega:\Omega\to G.$ As the domain $G$ one can take any plane
domain which is conformally equivalent to $\Omega.$

It is well known that $\omega\in W^{1,p}_{\rm loc}(\Omega),$ $2\leq
p<{2K\over K-1},$ $\omega\in C^\alpha_{\rm loc}(\Omega),$
$\alpha=1/K,$ and we have seen that if $\omega(z)=u(z)+i\, v(z),$
 then
 ${\rm div}\,(A(z)\nabla\,u(z))=0,$ ${\rm
div}\,(A(z)\nabla\,v(z))=0$ a.e. in $\Omega$ and $u$ and $v\in C\cap
W^{1,2}_{\rm loc}(\Omega).$ In what follows we will say that the
matrix function $A$ generates the corresponding quasiconformal
mapping $\omega,$ or that $A$ and $\omega$ {\it are agreed}.

Note also in passing the very often applied fact that quasiconformal
mappings $\omega$ admit the change of variables in integrals because
homeomorphisms of the class $W^{1,2}_{\rm loc}$ are absolute
continuous with respect to the area measure, see e.g. Theorem
III.6.1, Lemmas III.2.1 and III.3.3 in \cite{LV:book}.

\par
We complete this section with the following very important result on the composition operators in the Sobolev spaces.
\par
Let $\Omega$ be a domain in the $n-$dimensional Euclidean space
${\Bbb R}^n,$ $n\geq 2.$ For the sake of completeness in the
exposition we recall that the Sobolev space $L^1_ p(\Omega),$ $p\geq
1,$ consists of locally integrable functions $\varphi : \Omega\to
{\Bbb R}$  with distributional partial derivatives of the first
order with the seminorm
$$||\varphi||_{L_p^1(\Omega)}=||\nabla\varphi||_{L_p(\Omega)}=\left(\int\limits_\Omega
|\nabla\varphi|^p dm\right)^{1/p}<\infty,$$
where $m$ is the Lebesgue measure in ${\Bbb R}^n,$ $\nabla\varphi$ is the distributional gradient
of the function $\varphi,$
$\nabla\varphi=\left(\partial\varphi/\partial x_1,...,\partial\varphi/\partial x_n\right),$
$x=(x_1,...,x_n),$ defined by the
conditions
$$\int\limits_\Omega \varphi\cdot{\partial\eta\over\partial x_i}\,dm=
-\int\limits_\Omega {\partial\varphi\over\partial x_i}\cdot\eta\,dm
\,\,\,\,\forall\eta\in C_0^\infty(\Omega),\,\,i=1,2,...,n.$$
Here  $C_0^\infty(\Omega)$ stands for
 the space of all infinitely smooth functions with a
compact support in $\Omega.$ Similarly, a vector-function is said to belong to
the Sobolev class $L_p^1(\Omega)$
if its coordinate functions belong to $L_p^1(\Omega).$ The
classes $W^{1,p}(\Omega) = L_p^1(\omega)\cap L_p(\Omega),$
which will be used later on,
 differ from the classes $L_p^1(\Omega)$
only by
the norm
$$||\varphi||_{W^{1,p}(\Omega)} = ||\varphi||_{L_p(\Omega)}+||\nabla\varphi||_{L_p(\Omega)}.$$

\par
The following statement  will play an important role in our further
con\-si\-de\-ra\-tions, see \cite{Goldshtein}, \cite{U}  and
\cite{VU}.

\bigskip

{\bf Lemma 2.1.} {\it Let $\omega : \Omega\to G$ be a homeomorphism.
Then the following conditions are equivalent:
\par
1) the composition $\omega^*\varphi=\varphi\circ \omega$ generates
the bounded operator
$$\omega^*:\,L_p^1(G)\to L_q^1(\Omega)\ ,\,\,\,1\leq q\leq p<\infty\ ,$$

\par
2) the mapping $\omega$ belongs to the class $W^{1,1}_{\rm
loc}(\Omega)$ and the function
$$K_p(x; \omega) := \inf\left\{
k(x) : ||D\omega||(x)\leq k(x)|J_\omega (x)|^{1/p}\right\}$$ belongs
to $L_r(\Omega),$ where $r$ is determined from the equality
$${1\over r}={1\over q}-{1\over p}\ .$$}

Here $||D\omega||(x)$ stands for the operator norm of the Jacobian
matrix $D\omega$ of the mapping $\omega$ at the point $x,$ i.e.,
$||D\omega||(x) := \sup\limits_{h\in{\Bbb C},|h|=1}D\omega\cdot h.$
\par
Assuming that $\omega$ is a quasiconformal homeomorphism and,
con\-se\-qu\-ent\-ly, $\omega\in W^{1,2}_{\rm loc}(\Omega)$ and
$K_2(x,\omega)\in L^\infty(\Omega),$ we arrive at the the following
conclusion, when $n = 2,$ $p = 2,$ $q=2$ and $r =\infty.$

\bigskip

{\bf Proposition 2.1.} {\it Let $\omega:\Omega\to G$ be a
quasiconformal homeomorphism and let $\varphi: G\to {\Bbb C}$ belong
to the class $W^{1,2}_{\rm loc}(G).$ Then the composition function
$\varphi\circ \omega\in W^{1,2}_{\rm loc}(G). $}

\bigskip

The study of superposition operators on Sobolev spaces stems from
the classical article \cite{Sobolev}, see also e.g. \cite{GGR},
\cite{V} and \cite{VE} for the detailed history and bibliography.
The Reshetnyak problem on the description of all isomorphisms of the
Sobolev space $L^1_n$ generated by a quasiconformal mapping in
${\mathbb R}^n$, $n\ge 2,$ was posed in 1968 at the first Donetsk
colloquium in the theory of quasiconformal mappings. The problem was
solved in \cite{VG}, see also \cite{Reshetnyak}. As a consequence,
the composition $\omega^*\varphi=\varphi\circ \omega$ generates the
bounded operator
$$\omega^*:\,L_2^1(G)\to L_2^1(\Omega)$$
if and only if the homeomorphism $\omega$ is quasiconformal, see
also \cite{RR}.

\section{A basic identity}

It is well-known that every positive defined quadratic form
\begin{equation}\label{adm0} ds^2=a(x,y)dx^2+2b(x,y)dxdy+c(x,y)dy^2, 
\end{equation} defined
in a plane domain $\Omega,$ can be reduced, by means of a suitable
quasiconformal change of variables, to the canonical form
\begin{equation}\label{adm00} ds^2=\Lambda(du^2+dv^2),\,\, \Lambda\not=
0\,\,\,\mbox{a.e. in $\Omega,$} 
\end{equation} provided that
$ac-b^2\geq\triangle_0>0,$ $a>0,$ a.e. in $\Omega,$ see, e.g.,
\cite{BGMR}, pp. 10--12. This key result can be extended to every
linear divergent operator of the form ${\rm div\,}[A(z)\nabla
u(z)],$ $z=x+iy,$  for  matrix function $A\in M^{2\times
2}(\Omega).$ Namely, given such a matrix function $A$ and a
quasiconformal mapping $\omega:\Omega\to G,$ $\omega\in W^{1,2}_{\rm
loc}(\Omega)$ agreed with  $A,$  we have already seen by direct
computation, that if the function $T$ and the entries of  $A$ are
sufficiently smooth, then 
\begin{equation}\label{adm1} {\rm
div\,}[A(z)\nabla(T(\omega(z)))]=J_\omega(z)\triangle
T(\omega(z)),\,\,z\in\Omega, 
\end{equation} see, e.g., \cite{GNR-DAHY},
\cite{GNR2016}. Here $J_\omega(z)$ stands for the Jacobian of the
mapping $\omega(z),$ e.g., $J_\omega(z)={\rm det\,}D_\omega(z),$
where $D_\omega(z)$ is the Jacobian matrix  of the mapping $\omega$
at the point $z\in\Omega.$ Making use of the standard procedure, the
regularity requirements  in the equality
 (\ref{adm1}) may be substantially weakened. The equality (\ref{qk4}) below,
 that will be applied to the study of weak solutions to some semi-linear partial differential
equations, can be viewed as a weak counterpart to equality
(\ref{adm1}).

\bigskip

{\bf Proposition 3.1.} {\it Let $\Omega$ be a domain in ${\Bbb C}$,
$A\in M^{2\times 2}(\Omega)$ and $\omega:\Omega\to G$ be a
quasiconformal mapping agreed with  $A.$ Then the equality
\begin{equation}\label{qk4} \int_\Omega \langle
A(z)\nabla(T(\omega(z))),\,\nabla\varphi(z)\rangle\,dm_z
=\int_\Omega\langle D^{-1}_\omega(z)\nabla
T(\omega(z)),\,\nabla\varphi(z)\rangle J_\omega(z)\,dm_z 
\end{equation}
 holds
for  every $T\in W^{1,2}_{\rm loc}(G)$ and for all $\varphi\in
W^{1,2}_0(\Omega).$}

\bigskip

\begin{proof} Assuming that $T\in W^{1,2}_{\rm loc}(G)$  and that
$\omega:\Omega\to G$  is a quasiconformal mapping
 agreed with  $A(z),$ we see, by Proposition 2.1,
 that $u:\,=\,T\circ\omega\in W^{1,2}_{\rm loc}(\Omega).$
 Since
\begin{equation}
 \nabla u(z)=D^t_\omega(z)\nabla T(\omega(z)), 
 \end{equation} where
$D^t_\omega(z)$ stands  for the transpose matrix to $D_\omega(z),$
and $\omega$ satisfies the Beltrami equation (\ref{Beltrami2}),
that can be written in the matrix form as 
\begin{equation}
A(z)D^t_\omega(z)=D^{-1}_\omega(z)\,J_\omega(z), 
\end{equation} we arrive
successively to the required equality (\ref{qk4}):
\begin{multline}\label{qk1b}
\int_\Omega \langle
A(z)\nabla(T(\omega(z))),\,\nabla\varphi(z)\rangle\,dm_z\\=\int_\Omega\langle
A(z)D^t_\omega(z)\nabla T(\omega(z)),\,\nabla\varphi(z)\rangle\,dm_z
\\=\int_\Omega\langle D^{-1}_\omega(z)\nabla T(\omega(z)),\,\nabla\varphi(z)\rangle
J_\omega(z)\,dm_z.
\end{multline}
\end{proof}

\par
{\bf Remark 3.1.} It is easy to see that if $T$ and $A$ are smooth
and (\ref{qk4}) holds, then (\ref{adm1}) also  holds. One can use
the Green's
 formula as above, but in the opposite direction, to show that
\begin{equation}\label{adm11}
\int_\Omega\left\{{\rm
div\,}[A(z)\nabla(T(\omega(z)))]-J_\omega(z)\triangle
T(\omega(z))\right\}\varphi(z)\,dm_z=0,\,\,\forall\varphi\in
C^1_0(\Omega). 
\end{equation}
 Indeed, by Green's formula we have that
\begin{equation} \int_\Omega
\langle
A(z)\nabla(T(\omega(z))),\,\nabla\varphi(z)\rangle\,dm_z=-\int_\Omega
{\rm div\,}[A(z)\nabla(T(\omega(z)))]\varphi(z)\,dm_z.
\end{equation}
On the other hand,  applying again Green's formula after the change of variables $w=\omega(z),$ we get
\begin{multline}\label{qk1c}
\int_\Omega J_\omega(z)\triangle T(\omega(z))\varphi(z)\,dm_z=\\ -\int_\Omega\langle\nabla T(\omega(z)),
\,[D^t_\omega]^{-1}(z)\nabla\varphi(z)\rangle J_\omega(z)\,dm_z.
\end{multline}
From the linear algebra it follows that
\begin{equation}\label{mmm0}
\langle\nabla T(\omega(z)),\,J_\omega(z)[D^t_\omega]^{-1}(z)\nabla\varphi(z)\rangle=
\langle J_\omega(z)D^{-1}_\omega(z)\nabla T(\omega(z)),\,\nabla\varphi(z)\rangle.
\end{equation}
Indeed, let $\omega(z)=a(z)+ib(z),$ $z=x+iy,$ then
\begin{equation}\label{mmm1}
D_\omega=\left(\begin{array}{ccc} a_x  & a_y\\
                            b_x   & b_y
                              \end{array}\right),
\end{equation}
\begin{equation}\label{mmm2}
D^t_\omega=\left(\begin{array}{ccc} a_x  & b_x\\
                            a_y   & b_y
                              \end{array}\right),
                              \end{equation}
\begin{equation}\label{mmm3}
D^{-1}_\omega J_\omega=\left(\begin{array}{ccc} b_y  & -a_y \\
                            -b_x   & a_x
                              \end{array}\right)
                              \end{equation}
and
\begin{equation}\label{mmm4}
[D^t_\omega]^{-1}J_\omega=\left(\begin{array}{ccc} b_y  & -b_x \\
                            -a_y   & a_x
                              \end{array}\right).
                                                            \end{equation}
Setting $\nabla T=(u,v),$ and
$\nabla\varphi=(\xi,\eta),$   we  get that
\begin{multline}
\langle \nabla T,\,J_\omega[D^t_\omega]^{-1}\nabla\varphi\rangle=\langle(u,v),\,(b_y\xi-b_x\eta, -a_y\xi+a_x\eta)\rangle=\\
u(b_y\xi-b_x\eta)+v(a_x\eta-a_y\xi),
\end{multline}
\begin{multline}
\langle J_\omega D^{-1}_\omega\nabla T,\,\nabla\varphi\rangle=\langle(b_yu-a_yv, -b_xu+a_xv),\,(\xi,\eta)\rangle=\\
u(b_y\xi-b_x\eta)+v(a_x\eta-a_y\xi),
\end{multline}
and  we arrive to formula (\ref{mmm0}).

\par
Thus,
\begin{multline}\label{qk1ca}
\int_\Omega J_\omega(z)\triangle T(\omega(z))\varphi(z)\,dm_z=\\
-\int_\Omega\langle D^{-1}_\omega(z)\nabla T(\omega(z)),\,\nabla\varphi(z)\rangle J_\omega(z)\,dm_z.
\end{multline}

\par
Finally, since $C^1_0(\Omega)$ is dense in $W^{1,2}_0(\Omega),$ we
come to the equality (\ref{adm1}) for almost all points in $\Omega.$
However, since all functions are smooth, the equation holds at every
point.

\section{Semi-linear model elliptic equations}

Let $\Omega$ be a bounded domain in ${\Bbb C}$ and let $f:{\Bbb R}\to{\Bbb R}$ be a continuous function. In this section we
study a model semi-linear equation
\begin{equation}\label{3.1a}
 {\rm div\,}[A(z)\nabla u(z)]=f(u(z)), \,\,\,z\in\Omega,
\end{equation} as well as its Laplace's counterpart: 
\begin{equation}\label{3.3a}
\triangle T(w)=J(w)f(T(w)),\,\,\,w\in G=\omega(\Omega),
\end{equation} where $\omega:\Omega\to
G$ is a quasiconformal mapping
 agreed with  $A(z)$ and
$J(w)$ stands for the Jacobian of the inverse mapping
$\omega^{-1}:G\to\Omega.$
\par
Under {\it a weak solution to the equation} (\ref{3.1a})  we
understand a function $u\in C\cap W^{1,2}_{\rm loc}(\Omega)$ such
that 
\begin{equation}
\label{qk10} \int_\Omega \langle A(z)\nabla
u(z),\,\nabla\varphi(z)\rangle\,dm_z+\int_\Omega
f(u(z))\varphi(z)\,dm_z=0\,\,\,\forall \varphi\in C\cap
W^{1,2}_0(\Omega)\ . 
\end{equation}
  \par
 The equation (\ref{3.3a}) contains a weight function $J(w)\in L^1_{\rm loc}(G).$ In spite of this, we may define a notion of a weak solution
 to the equation (\ref{3.3a}) in the similar way.
\par
 We say that $T$ is {\it a weak solution to the equation} (\ref{3.3a}) if $T\in C\cap W^{1,2}_{\rm loc}(G)$
 and
\begin{equation}\label{qk11} \int_G\langle\nabla T(w),\,\nabla
\psi(w)\rangle\,dm_w+\int_G J(w)f(T(w))\psi(w)\,dm_w=0\,\,\,\forall
\psi\in C\cap W^{1,2}_0(G)\ . 
\end{equation} Since $J(w)$ is the Jacobian of
the mapping $\omega^{-1}(w),$ it is easy to verify, by performing
the change of variable by the formula $w=\omega(z),$ that the second
integral in (\ref{qk11}) is well-defined. Here  we again made use of
the fact that the composed mapping $u(z)=T(\omega(z))$ is  in $C\cap
W^{1,2}_{\rm loc}(\Omega)$ if  $T\in C\cap W^{1,2}_{\rm loc}(G)$ and
$\omega$ is quasiconformal.
\par
The existence of weak solutions to the equations (\ref{3.1a}) and
(\ref{3.3a} ), under specified conditions on the right hand sides, as well as fundamental properties of solutions can be found, e.g. in \cite{LU},
see also \cite{Skrypnik}.

\bigskip

The following Factorization Theorem is a principal result of this
paper.

\bigskip

{\bf Theorem 4.1.} {\it Let $\Omega$ be a domain in ${\Bbb C}$,
$A\in M^{2\times 2}(\Omega)$ and  let $f:{\Bbb R}\to{\Bbb R}$ be a
continuous function. Then every weak solution $u$ of the semi-linear
equation 
\begin{equation}\label{qk12a}
 {\rm div\,}[A(z)\, \nabla u(z)]\ =\ f(u(z)), \,\,\,z\in\Omega,
\end{equation} can be represented as the composition 
\begin{equation}\label{factor}
 u(z)=T(\omega(z)), 
 \end{equation} where
$\omega:\Omega\to G$ is a quasiconformal mapping agreed with $A$ and
$T$ is a weak solution to the equation 
\begin{equation}\label{eq4a}
\triangle\,T(w)\ =\ J(w)\, f(T(w)),\,\,\, w\in G\ . 
\end{equation} Here $J(w)$
stands for the Jacobian of the inverse mapping $\omega^{-1}(w).$ }

\bigskip

\begin{proof} Let $u$ be a weak solution  of the semi-linear equation
(\ref{qk12a}) and $T=u\circ\omega^{-1}$. Then by Proposition 2.1
$T\in C\cap W^{1,2}_{\rm loc}(G)$ and we have that 
\begin{equation}\label{qk14}
\int_\Omega \langle A(z)\nabla
(T(\omega(z))),\,\nabla\varphi(z)\rangle\,dm_z+\int_\Omega
f(T(\omega(z)))\varphi(z)\,dm_z=0 
\end{equation} for all $\varphi\in C\cap
W^{1,2}_0(\Omega).$ Next, by Proposition 3.1, 
\begin{equation}\label{qk15}
\int_\Omega \langle
A(z)\nabla(T(\omega(z))),\,\nabla\varphi(z)\rangle\,dm_z=
\int_\Omega\langle D^{-1}_\omega(z)\nabla
T(\omega(z)),\,\nabla\varphi(z)\rangle J_\omega(z)\,dm_z, 
\end{equation} and
therefore 
\begin{equation}\label{qk15a} \int_\Omega\langle J_\omega(z)
D^{-1}_\omega(z)\nabla T(\omega(z)),\,\nabla\varphi(z)\rangle
\,dm_z+ \int_\Omega f(T(\omega(z)))\varphi(z)\,dm_z=0 
\end{equation} for all
$\varphi\in C\cap W^{1,2}_0(\Omega).$

\par
Given an arbitrary function $\psi(w)\in C\cap W^{1,2}_0(G),$ we can
set in (\ref{qk14}) and (\ref{qk15}) $\varphi(z)=\psi(\omega(z)),$
because, by Proposition 2.1,
 such $\varphi\in C\cap W^{1,2}_0(\Omega).$
Performing in (\ref{qk15a})
the change of variable
by the formula $z=\omega^{-1}(w),$ we obtain
\begin{multline*}\label{qk17}
\int_G\langle
J_\omega(\omega^{-1}(w))D^{-1}_\omega(\omega^{-1}(w))\nabla
T(w),\,D^t_\omega(\omega^{-1}(w))\nabla\psi(w)\rangle J(w) \,dm_w\\+
\int_G J(w)f(T(w))\psi(w)\,dm_w=0.
\end{multline*}
 Since, by elementary algebraic arguments,
$$
\langle J_\omega(\omega^{-1}(w))D^{-1}_\omega(\omega^{-1}(w))\nabla
T(w),\,D^t_\omega(\omega^{-1}(w))\nabla\psi(w)\rangle =
$$
$$
=J_\omega(\omega^{-1}(w))\langle\nabla T(w),\,\nabla\psi(w)\rangle,
$$ and
$$J_\omega(\omega^{-1}(w))=1/J(w),$$
we see that the identity 
\begin{equation}\label{qka11} \int_G\langle\nabla
T(w),\,\nabla \psi(w)\rangle\,dm_w+\int_G J(w)f(T(w))\psi(w)\,dm_w=0
\end{equation} holds for all $\psi\in C\cap W^{1,2}_0(G).$  Thus,  $T$ is a
weak solution to the equation (\ref{eq4a}).
\end{proof}

\par
{\bf Remark 4.1.}\, Since the arguments given above are invertible,
we see that if $T$ is a weak solution to the equation (\ref{eq4a}),
then the function $u(z)=T(\omega(z))$ is a weak solution to the
semi-linear equation (\ref{qk12a})

\bigskip

Assuming that the function $f$ is non-negative, we arrive at the following statement.

\bigskip

{\bf Corollary 4.1.} {\it Let $\Omega$ be a domain in ${\Bbb C}$,
$A\in M^{2\times 2}(\Omega)$ and  let $f:{\Bbb R}\to{\Bbb R}$ be a
continuous function. If $f(u)\geq 0,$ then  every weak solution $u$
of the semi-linear equation 
\begin{equation}\label{qk12ab}
 {\rm div\,}[A(z)\nabla u(z)]=f(u(z)), \,\,\,z\in\Omega,
\end{equation} can be represented as the composition \begin{equation}\label{factor}
u(z)=T(\omega(z)), \end{equation} where $\omega:\Omega\to G$ is quasiconformal
mapping agreed with $A$ and $T$ is a subharmonic function, being a
weak solution of the equation \begin{equation}
\label{eq4ab} \triangle\,T\ =\
J(w)\,f(T(w))\ \ \ \mbox{ in\ $G$}\ . 
\end{equation} Here $J(w)$ stands for the
Jacobian of the inverse mapping $\omega^{-1}(w).$ }

\bigskip

Among the quasiconformal mappings $\omega:\Omega\to G$ there are a
variety of the so-called volume-preserving maps, for which
$J_\omega(z)\equiv 1,$ $z\in \Omega.$  Many examples of such maps
will be given in the next section. The following statement  may have
of independent interest.

\bigskip

{\bf Corollary 4.2.} {\it Let $A\in M^{2\times 2}(\Omega)$ be a
matrix function that generates a volume-preserving quasiconformal
mapping $\omega :\Omega\to G.$ Then every weak solution $u$ of the
semi-linear equation \begin{equation}\label{qk12}
 {\rm div\,}[A(z)\nabla u(z)]=f(u(z)), \,\,\,z\in\Omega,
\end{equation} can be represented as the composition \begin{equation}\label{factor}
u(z)=T(\omega(z)), \end{equation} where $T$ is a weak solution to the equation
\begin{equation}\label{eq4} \triangle\,T=f(T(w)),\,\,\mbox{a.e. in $G$}. \end{equation}  }

\bigskip

{\bf Remark 4.2.} By the Measurable Riemann mapping theorem, see
e.g. Theorem V.B.3 in \cite{Ahlfors:book} and Theorem V.1.3 in
\cite{LV:book}, given $\mu(z),$ $z\in\Omega,$ agreed with the matrix
function $A\in M^{2\times 2}(\Omega),$ there exists a quasiconformal
mapping $\omega:\Omega\to G$ with the complex dilatation $\mu$ given
by (2.6). Here if $\Omega$ is finitely connected, then  $G$ is
either the canonical circular domain, or the complex plane, see e.g.
Theorem V.6.2 in \cite{Golusin}. It allows us to choose in the
Factorization Theorem 4.1 as the domain $G$ one from the above
canonical domains. If $\Omega$ is simply connected with a
nondegenerate boundary then we can set $G={\Bbb D}.$ The latter
makes it possible to remove the restrictions on the regularity of
the boundary in the study of boundary value problems for equation
(\ref{eq4a}), see Section 8.

\bigskip

\section{Some   Applications}

We start with a few examples of the application of Theorem 4.1 to
the study  of boundary blow-up solutions for some classical model
semi-linear elliptic equations.
\par
Let $\Omega$ be a bounded domain in ${\Bbb C}$ and let
$\partial\Omega$ denote its boundary. In this section we study the
problem \begin{equation}\label{3.1} {\rm div\,}[A(z)\nabla u(z)]=f(u(z)), 
\end{equation}
\begin{equation}\label{3.2} u(z)\to\infty\ \ \ \  \mbox{as\ \ \ \ $d(z):={\rm
dist}\,(z,\partial\Omega)\to 0$}, 
\end{equation} as well as its Laplace
counterpart: \begin{equation}\label{3.3} \triangle u(z)=f(u(z)), \end{equation}
\begin{equation}\label{3.4} u(z)\to\infty\ \ \ \  \mbox{as\ \ \ \ $d(z):={\rm
dist}\,(z,\partial\Omega)\to 0$}. \end{equation} Solutions to these problems
are called {\it boundary blow-up solutions,} or {\it large
solution.}
\par
The existence of a large solution to (\ref{3.3}) is related to the
existence of a maximal solution $\tilde u$ of (\ref{3.3}) in
$\Omega,$ which in turn depends on the so called Keller-Osserman
condition, see \cite{Keller} and \cite{Oss}. For example, in the
paper \cite{Keller} a simple upper bound was obtained for any
solution, in any number of variables, of semi-linear equation
(\ref{3.3}). The bound depends on the function $f$ which must be
positive and satisfy the well-known Keler-Osserman condition. Recall
that a function $f\in C({\Bbb R}_+)$ satisfies the Keler-Osserman
condition if there exists a positive non-decreasing function $h$
such that
\begin{equation}\label{ko} f(t)\geq h(t), \forall t\in {\Bbb R}_+
\,\,\mbox{and}\,\,\,\int_{t_0}^\infty\left\{\int_0^t
h(s)ds\right\}^{-1/2}\,dt<\infty \,\,\mbox{for all $t_0>0.$} \end{equation} It
is known that if $f$ is non-decreasing and satisfies the
Keller-Osserman condition, then a {\it large solution exists in
every bounded smooth domain.} Uniqueness in smooth domains was
established under some additional conditions on $f,$ see e.g.
\cite{Bandle_Marcus1995}. It is easy to check that the functions
$f(t)=e^t$ and $f(t)=t^p,$ $p>1,$ satisfy (\ref{ko}). The
Gauss-Bieberbach-Rademacher semi-linear equation \begin{equation}\label{LB}
\triangle u=e^u, \end{equation} as far as we know,  was first investigated by
Bieberbach in his pioneering work \cite{Bieberbach} related to the
study of automorphic functions in the plane. It is this work has
stimulated numerous studies in the field of semi-linear differential
equations in $R^n,$ $n\geq 1,$ and the equation (\ref{LB}) continues
to play the role of one of the fundamental model equations of the
theory. It is important to note that in simply connected planar
domains $\Omega$ the large solutions for the equation (\ref{LB}) are
represented explicitly by means of the conformal map $\omega :
\Omega\to \mathbb D:=\{ z\in\mathbb C: |z| < 1\}:$ \begin{equation}\label{L} u(z)
= \log{8|\omega'(z)|^2\over (1-|\omega(z)|^2)^2}. \end{equation} Let us also
remark that, given an arbitrary analytic function $\omega(z)$ in
$\Omega,$ the formula (\ref{L}) generates a solution, generally
speaking with singularities, of the equation (\ref{LB}).
\par
For the model case of
equation
\begin{equation}\label{1.5} \Delta u=e^{au},\,\,\, a>0,
\end{equation}
the
following result holds, see \cite{ML2014}, Theorem 5.3.7.

\bigskip

{\bf Theorem 5.1}\,\, {\it Let $\Omega$ be a bounded
domain in ${\Bbb C}$ such that
$\Omega={\rm Int}\,{\overline\Omega}.$ Then there exists
one and only one blow-up solution to (\ref{1.5}).}

\bigskip

The Gauss-Bieberbach-Rademacher semi-linear equation is one of the
principal model equations in the theory of non-linear partial
differential
 equations and their applications, see e.g. \cite{Neklyudov}, \cite{ML2014}, \cite{Oleynik} and the references therein.
 Note that the equation
appears also as a model one in problems of differential geometry in relation with existence of surfaces
of negative Gaussian curvature \cite{Vekua} and in studying the
equilibrium of a charged gas.
\par
In the general case for $|\mu(z)|\leq k<1$ the solution of the
Beltrami equation can be written as an infinite series of singular
integral transforms of Hilbert and Cauchy type for complex
dilatations, see, e.g. \cite{BGMR}, p. 33. Here we give the
explicit solutions of the Beltrami equation for the cases where
the complex dilatation is a measurable function that depends on a
single real variable $x={\rm Re}\, z$, $y={\rm Im}\,z$, $\arg z$
or $|z|,$ see \cite{BGMR}, \S 5.10, \cite{GR1995}.
 We make use of the corresponding explicit formulae
to write down a number of explicit  solutions for the counterpart
(5.1) of the Gauss-Bieberbach-Rademacher equation in the unit disk
and the upper half-plane.
\par
 We start with the following statement, see,
e.g., \cite{BGMR}, p. 82.

\bigskip

{\bf Proposition 5.2.} {\it  Let the complex dilatation $\mu(z)$ has
 the form
\begin{equation}\label{eq1.50}
 \mu(z)=k(|z|){z\over{\bar z}}\,,
 \end{equation}
where $k(\tau):{\Bbb R}\to{\Bbb C}$ is a measurable function, $\|
k\|_{\infty}<1$. Then the formula \begin{equation}\label{eq1.51}
\omega(z)={z\over{|z|}}\exp
\left\{-\int\limits_{|z|}^{1}{{1+k(\tau)}
\over{1-k(\tau)}}\,{{d\tau}\over{\tau}}\right\}\end{equation} represents a
unique quasiconformal mapping of the unit disk onto itself with the
complex dilatation $\mu$ and the normalizations  $\omega(0)=0$ and
$\omega(1)=1.$ }

\bigskip

Analyzing formula (\ref{eq1.51}), we arive at the following
statement.

\bigskip

{\bf Corollary 5.1.} {\it Let ${\Bbb D}$ be  the unit disk in the complex
plane ${\Bbb C}$ centered at the origin and
 let the matrix function $A(z)$ be generated by  the Beltrami coefficient
\begin{equation}\label{mu1} \mu(z)=k(|z|){z\over\bar z} \end{equation} in accordance with
formula (\ref{matrix}) where \begin{equation}\label{vp} k(t)=\nu^2(t)\pm
i\,\nu(t)\sqrt{1-\nu^2(t)} \end{equation} and $\nu(t),$ $0\leq t<1,$ stands
for an arbitrary measurable real-valued function such that
$|\nu(t)|\leq q<1.$ Then there exists one and only one boundary
blow-up solution to the  semi-linear equation \begin{equation}\label{eq1.50a}
{\rm div\,}[A(z)\nabla u]=e^{u},\,\,\,z\in {\Bbb D,}
\end{equation} which is written
explicitly by the Liouville-Bieberbach formula \begin{equation}\label{solution}
u(z)=\log{8\over (1-|z|^2)^2}. \end{equation} }
\par\medskip

\begin{proof} Let $A(z)$ be a matrix function generated by the complex
Beltrami coefficient (\ref{mu1}), satisfying (\ref{vp}). By the
Measurable Riemann mapping theorem, see e.g. Theorem V.B.3 in
\cite{Ahlfors:book} and Theorem V.1.3 in \cite{LV:book}, there
exists the unique normalized quasiconformal self-homeomorphism
$\omega$ of the unit disk $\mathbb D$ with the Beltrami coefficient
$\mu(z)=k(|z|)z/\bar z$ and this homeomorphism, by Proposition 5.2,
is written explicitly, \begin{equation}\label{explic} \omega(z)={z\over
|z|}\exp\left\{-\int_{|z|}^1{1+k(t)\over 1-k(t)}{dt\over t}\right\}.
\end{equation} We see that
$${\omega_{\bar z}\over \omega}={k\over 1-k}{1\over \bar z},$$
$${\omega_z\over \omega}={1\over 1-k}{1\over z}$$
and therefore, for the Jacobian
$J_\omega(z)=|\omega_z|^2-|\omega_{\bar z}|^2$ we  have
\begin{equation}\label{Jac}
J_\omega(z)={1-|k|^2\over|1-k|^2}{|\omega|^2\over|z|^2}. \end{equation} In
order to apply Theorem 4.2, we have to verify that
$J_\omega(z)\equiv 1$ and $|\omega(z)|=|z|$ for $z\in {\Bbb D}.$ Indeed,
noting that ${\rm Re}\,k=|k|^2$ and substituting $k(t)$ given by
formula (\ref{vp}) into (\ref{explic}) we see that
$|\omega(z)|=|z|,$ $z\in {\Bbb D}.$ The later implies that
$J_\omega(z)\equiv 1.$
\par
 Making use of
Corollary 4.2, we see that $u(z)=T(\omega(z))$  and $T$ satisfies
the Gauss-Bieberbach-Rademacher equation \begin{equation}\label{semi-linear}
\triangle T=e^T\,\,\, \mbox{in ${\Bbb D}.$} \end{equation} By the mentioned
above Bieberbach's result  and referring also to  Theorem 5.3.7 from
\cite{ML2014}, we see that the unique  boundary blow-up solution to
the equation (\ref{semi-linear}) in the unit disk ${\Bbb D}$ is
given by the explicit formula \begin{equation} T(w)=\log{8\over(1-|w|^2)^2}. \end{equation}
Since $|\omega(z)|=|z|$ for $z\in {\Bbb D},$ the required solution
to the equation (\ref{eq1.50a}) also has the explicit representation
\begin{equation} u(z)=\log{8\over (1-|z|^2)^2}. \end{equation}
\end{proof}

\par
Thus, there is an infinite set of matrix functions $A(z)$ such that
the corresponding anisotropic semi-linear equations of the form
(\ref{eq1.50a}) have a solution defined by the explicit formula
(\ref{solution}).
\par
Let us give a nontrivial example, based on the well-known
logarithmic spiral quasiconformal mapping
$$\omega(z)=ze^{2i\log|z|}$$
which plays important role in the study of different problems of
contemporary analysis, see, e.g., \cite{BGMR}, \S 13.2, \cite{GM},
\cite{GM:2001}. This function $\omega$ maps the unit disk ${\Bbb D}$
onto itself and transforms radial lines into spirals, infinitely
winding about the origin, and it is just the volume preserving. The
mapping $\omega$ satisfies the Beltrami equation with
$$\mu(z)={\omega_{\bar z}\over \omega_z}={1\over 2}(1+i){z\over\bar z}$$
and  the Jacobian $J_{\omega}(z)=|\omega_{z}(z)|^2-|\omega_{\bar
z}(z)|\equiv 1.$ We see that $\mu$ corresponds to (\ref{vp}) with
$\nu(t)\equiv 1/\sqrt 2.$ Since
$${\rm Re}\,\mu(z)={x^2-y^2-2xy\over 2(x^2+y^2)},\,\,\,{\rm Im}\,\mu(z)={x^2-y^2+2xy\over 2(x^2+y^2)},\,\,z=x+iy,$$
we see that $\mu$ generates by the formula (\ref{matrix}) the
matrix function $A_{sp}(x,y)$ with the following entries:
$$a_{11}=\alpha={|1-\mu|^2\over 1-|\mu|^2}=3-2{x^2-y^2-2xy\over (x^2+y^2)},$$
$$a_{22}=\delta={|1+\mu|^2\over 1-|\mu|^2}=3+2{x^2-y^2-2xy\over (x^2+y^2)},$$
$$a_{12}=a_{21}=\beta={-2{\rm Im}\,\mu\over 1-|\mu|^2}=-2{x^2-y^2+2xy\over (x^2+y^2)}.$$
By Corollary 5.1 the semi-linear equation \begin{equation}\label{eq1.50b} {\rm
div\,}[A_{sp}(x,y)\nabla u]=e^{u},\,\,\,z\in {\Bbb D}, \end{equation} must have the
function
$$u(x,y)=-2\log(1-x^2-y^2)+\log 8$$
as the blow-up solution in the disk ${\Bbb D}.$ Let us verify this
conclusion by means of direct computation.
\par
We see that
$$u_x={4x\over 1-x^2-y^2},\,\,u_y={4y\over 1-x^2-y^2}$$
and therefore
$$\alpha u_x+\beta u_y=4{x+2y\over 1-x^2-y^2},\,\,\,\beta u_x+\delta u_y=4{y-2x\over 1-x^2-y^2}.$$
Then
$${\rm div}\,\left(A_{sp}(x,y)\nabla u\right)=4(\alpha u_x+\beta u_y)_x^{'}+4(\beta u_x+\delta u_y)_y^{'}={8\over(1-x^2-y^2)^2}=e^u,$$
and thus we successfully complete the verification.
\par
Note also that  the matrix function $A_{sp}$ in the polar coordinates
$z=re^{i\varphi}$ has the form \begin{equation}\label{matrix3}
A=\left(\begin{array}{ccc} 3-2\sqrt 2\cos(2\varphi+\pi/4)  & -2\sqrt 2\sin(2\varphi+\pi/4) \\
                           -2\sqrt 2\sin(2\varphi+\pi/4)    &  3+2\sqrt 2\cos(2\varphi+\pi/4)
                              \end{array}\right).
\end{equation}

\bigskip

{\bf Corollary 5.2.} {\it Let $\Omega$ be  the annulus $r<|z|<1$
in the complex plane ${\Bbb C}$ and
 let the matrix function $A(z)$ be generated by  the Beltrami coefficient
\begin{equation} \mu(z)=k(|z|){z\over\bar z} \end{equation} in accordance with  formula
(\ref{matrix}) where \begin{equation}\label{k2} k(t)=\nu^2(t)\pm
i\,\nu(t)\sqrt{1-\nu^2(t)} \end{equation} and $\nu(t),$ $0\leq t<1$ stands for
an arbitrary measurable function. If $|\nu(t)|\leq q<1,$ then there
exists one and only one boundary blow-up solution to the classical
model semi-linear equation \begin{equation} {\rm div\,}[A(z)\nabla
u]=e^{u},\,\,\,\mbox{in the annulus $r<|z|<1,$} \end{equation} and which is
written explicitly by the formula \begin{equation}\label{solution1}
u(z)=\log{2\pi^2\over|z|^2(\log^2r)\cdot\sin^2({{\pi\over\log
r}}\log|z|).} \end{equation}
}

\bigskip

\begin{proof} Arguing in the same way as in the proof of Corollary 5.1, we
see that $u(z)=T(\omega(z)),$ $r<|z|<1,$ where quasiconformal
self-homeomorphism $\omega(z)$ of the annulus is given by
(\ref{explic}), (note that under the condition (\ref{k2})
$|\omega(z)|=|z|$), and $T$ is the unique boundary blow-up solution
to the semi-linear equation \begin{equation}\label{semi-linear1} \triangle T=e^T
\end{equation} in the annulus $r<|w|<1.$ By Bieberbach's result, the unique
boundary blow-up solution to the equation (\ref{semi-linear1}) is
given by the formula \begin{equation}
T(\omega)=\log{8|F'(\omega)|^2\over(1-|F(\omega)|^2)^2} \end{equation} where
$F(\omega)$ stands for a conformal mapping of the annulus
$r<|\omega|<1$ onto the unit disk. It remains to find the
corresponding conformal mapping $F.$ We see that: 1) $w=\log\omega$
maps the annulus onto the strip $\log r<{\rm Re}\,w<0;$ 2)
$\zeta=-ie^{{\pi\over \log r}iw}$ maps the strip onto the right
half-plane; 3) $t=(\zeta-1)/(\zeta+1)$ maps the right half-plane
onto the unit disk. Composing  the above mappings we get the
required formula \begin{equation} F(\omega)={ie^{{\pi\over \log
r}i\log\omega}+1\over ie^{{\pi\over \log r}i\log\omega}-1}. \end{equation}
Next, if we set $\tau(\omega)=ie^{{\pi\over \log r}i\log\omega},$
then we see that \begin{equation} (1-|F(\omega)|^2)^2={16{\rm
Re}^2\,\tau(\omega)\over|\tau-1|^4} \end{equation} where
$${\rm Re}\,\tau(\omega)=-|\tau|\sin\left({{\pi\over \log r}\log|\omega|}\right).$$
On the other hand,
\begin{equation}
|F'(\omega)|^2={4\over|\tau-1|^4}|\tau'(\omega)|^2={4\over|\tau-1|^4}|\tau|^2{\pi^2\over
\log^2r}{1\over|\omega|^2}.
\end{equation}
Thus
\begin{equation}\label{annulus}
T(\omega)=\log{2\pi^2\over|\omega|^2(\log^2r)\cdot\sin^2({{\pi\over\log
r}}\log|\omega|)}
\end{equation}
and since $|\omega(z)|=|z|,$ we complete the
proof of Corollary 5.2.
\end{proof}

\par
Making use of the limit in (\ref{solution1}) as $r\to 0$ we get
the following result which may have of independent interest.

\bigskip

{\bf Corollary 5.3} {\it The semi-linear equation \begin{equation} {\rm
div\,}[A(z)\nabla u]=e^{u} \end{equation} for each matrix function $A$ from
Corollary 5.2 as well as the Gauss-Bieberbach-Rademacher equation
\begin{equation} \triangle u=e^u \end{equation} admit the following boundary blow-up
solution in the punctured unit disk $0<|z|<1$ \begin{equation}\label{news}
u(z)=\log{2\over |z|^2\log^2|z|}. \end{equation} }

\par\medskip
{\bf Remark 5.1.} The approach given above to the construction of a
boundary blow-up solution to the Gauss-Bieberbach-Rademacher
equation in the unit disk ${\Bbb D}$ with a singularity at the
origin can be extended to the case of a finite number of singular
points $z_k,$ $|z_k|<1,$ $k=1,2,...,n.$ Indeed, let $r>0$ be such
that all the disks $d_k=\{z:|z-z_k|\leq r\}$ belong to ${\Bbb D}$
and do not intersect each other. Denote by $F_r(z)$ a conformal
mapping of the circular multi-connected domain ${\Bbb
D}\setminus\cup_{k=1}^n d_k$ onto the unit disc ${\Bbb D},$ see
Theorem VI.1 in \cite{Golusin} due to Poincare. Then the required
solution with prescribed singularities at the points $z_k$ is given
by \begin{equation}\label{solutionS} u(z)=\lim_{r\to 0}\log{8|F'_r(z)|^2\over
(1-|F_r(z)|^2)^2}. \end{equation}
\par
In order to give more examples of applications of our Factorization
Theorem and its corollaries, recall the following statement, which
provides us with the explicit solutions of the Beltrami equation for
the cases where the complex dilatation $\mu$ is a measurable
function that depends on a single real variable $x={\rm Re}\, z$ or
$y={\rm Im}\,z,$ see, \cite{BGMR}, p. 78, \cite{GR1995}.

\bigskip

{\bf Proposition 5.3.} {\it Let $\mu:{\Bbb C}\to {\Bbb D}$ be an
arbitrary measurable function with $\Vert\mu\Vert_{\infty}\leq q<1$
that depends on $x={\rm Re}\,z$ only and let \begin{equation}\label{eq1.11}
\varphi(x)=\int\limits_{0}^{x}{{1+\mu (t)}\over{1-\mu(t)}}\,dt\,.
\end{equation} Then the formula \begin{equation}\label{eq1.12} \omega(z)=\varphi(x)+iy \end{equation}
represents a unique quasiconformal mapping of $\overline{\Bbb C}$
onto itself with complex dilatation $\mu$ and normalizations:
\begin{equation}\label{eq1.13}\omega(0)=0,\quad \omega(i)=i,\quad
\omega(\infty)=\infty\,. \end{equation}}

Indeed, from formulas (\ref{eq1.11}) and (\ref{eq1.12}) it follows
that
\begin{equation}\label{eq1.17}K^{-1}|x_1-x_2|\leq|\omega(z_1)-\omega(z_2)|\leq
K|z_1-z_2|\,. \end{equation} Thus, $\omega$ is a lipschitz homeomorphism of the
plane and hence ACL. The normalization (\ref{eq1.13}) is obvious.
Furthermore,
\begin{equation}\label{eq1.18}\omega_x=\varphi'(x)=\frac{1+\mu(x)}{1-\mu(x)}\,,
\end{equation}
\begin{equation}\label{eq1.19}\omega_y=i\,.\end{equation} Therefore,
\begin{equation}\label{eq1.20}\omega_{\bar z}={1\over
2}(\omega_x+i\omega_y)={\mu(x)\over{1-\mu(x)}}\,,\end{equation}
\begin{equation}\label{eq1.21}\omega_z={1\over2}(\omega_x-i\omega_y)={1\over{1-\mu(x)}}\,,\end{equation}
and hence $\omega$ satisfies the Beltrami equation
\begin{equation}\label{eq1.22}\omega_{\bar z}=\mu(x)\omega_z\,.\end{equation} The Jacobian
of $\omega$
\begin{equation}\label{eq1.23}J_\omega(z)={{1-|\mu (x)|^2}\over{|1-\mu
(x)|^2}} \geq K^{-1}>0
\end{equation}
 is positive, that is, $\omega$ is
orientation-preserving. This completes the proof of Proposition 5.3.

\bigskip

{\bf Remark 5.2.}  From \begin{equation}\label{eq1.14}u={\rm Re}\, \omega(z)={\rm
Re}\,\varphi(x)\,, \end{equation} \begin{equation}\label{eq1.15} v={\rm Im}\,
\omega(z)=y+{\rm Im}\,\varphi(x)\,, \end{equation} we conclude that  $\omega$
maps vertical lines onto vertical lines without contractions or
dilations, and the imaginary axis is mapped identically onto itself,
since $\varphi(0)=0$. It is easy to show that these geometric
properties characterize the class of quasiconformal mappings in
question, because
\begin{equation}\label{eq1.16}\mu(z)=\frac{\varphi'(x)-1}{\varphi'(x)+1}\end{equation}
depends only on $x.$
\par

\medskip

{\bf Corollary 5.4.} {\it Suppose that the complex dilatation $\nu
(z)$ of a quasiconformal mapping $g(z)$ of the extended complex
plane, keeping the points $0,1,\infty,$ depends only on $y = {\rm
Im}\,z$. Then \begin{equation}\label{eq1.24} g(z)=x+i\psi(y)\,, \end{equation} where
\begin{equation}\label{eq1.25}\psi(y)=\int\limits_{0}^{y} {{1-\nu(it)}\over
{1+\nu(it)}}\,dt\,.\end{equation}}

\bigskip

Indeed, let \begin{equation}\label{eq1.26} \omega={ \mathcal{A}}\circ g
\circ{ \mathcal{A}}^{-1}\,, \end{equation} where \begin{equation}\label{eq1.27}
{  \mathcal{A}}(\zeta)=e^{i{\pi\over2}}\zeta=i\zeta
\end{equation} is a counter-clockwise
rotation by the angle $\pi /2$. Then $\omega(0)=0$, $\omega(i)=i$,
$\omega(\infty)=\infty$, and the complex dilatation of $\omega$,
\begin{equation}\label{eq1.28}\mu(z)=-\nu(-iz)\,, \end{equation} depends only on $x={\rm
Re}\, z$. Thus, from (\ref{eq1.11}), (\ref{eq1.12}), and
(\ref{eq1.26}) we obtain (\ref{eq1.24}) and (\ref{eq1.25}).

\bigskip

Let us return to the study of  the blow-up solutions to the
Liouville-Biberbach type equation defined in the right half-plane.

\bigskip

{\bf Corollary 5.5.} {\it Let $H^+$ be  the right  half-plane
$\{z:{\rm Re}\, z>0\}$ in the complex plane ${\Bbb C}$  and
 let the matrix function $A(z)$ be generated by the formula (\ref{matrix}) with the Beltrami coefficient
\begin{equation}\label{eq1.30} \mu(x)=\nu^2(x)\pm i\,\nu(x)\sqrt{1-\nu^2(x)} \end{equation}
where $\nu(x),$ $x\in{\Bbb R},$ stands for an arbitrary measurable
real-valued function, such that $|\nu(x)|\leq q<1.$ Then there exist
boundary blow-up solutions to the  semi-linear equation
\begin{equation}\label{eq1.30a} {\rm div\,}[A(z)\nabla u]=e^{u},\,\,\,z\in H^+,
\end{equation} and which are written explicitly: \begin{equation}\label{eq1.31}
u(z)=\log{2\over x^2},\,\,z=x+iy, \end{equation} \begin{equation}\label{eq1.32}
u(z)=\log8\lambda^2-2\lambda x-2\log(1-e^{-2\lambda
x}),\,\,\lambda>0. \end{equation} }
\par\medskip

\begin{proof} We prove only that the function $u(z),$ given by the
formula (\ref{eq1.31}) solves the equation (\ref{eq1.30a}). Let
$A(z)$ be a matrix function generated by the complex Beltrami coefficient
$\mu,$ satisfying (\ref{eq1.30}). Since $\mu$ depends on $x$ only,
then formulae (\ref{eq1.11}), (\ref{eq1.12}) represent a unique
quasiconformal mapping $w=\omega(z)$ of the right half-plane $H^+$
onto itself with complex dilatation $\mu$ normalized by
$\omega(0)=0,$ $\omega(i)=i,$ $\omega(\infty)=\infty.$
\par
  In order to apply
Corollary 4.2, we have first to verify that $J_\omega(z)\equiv 1.$
Indeed, in our case
 ${\rm Re}\,\mu=|\mu|^2$ and formula (\ref{eq1.23})
implies that $J_\omega(z)\equiv 1.$ Moreover, ${\rm
Re}\,\omega(z)={\rm Re}\,\varphi(x)=x.$ By Corollary 4.2 a weak
solution $u$ of the semi-linear equation \begin{equation}\label{qk12ab}
 {\rm div\,}[A(z)\nabla u(z)]=e^u, \,\,\,z\in H^+,
\end{equation} is represented as the composition \begin{equation}\label{factor5}
u(z)=T(\omega(z)), \end{equation} where $T$ is a weak solution to the equation
\begin{equation}\label{eq4ab} \triangle\,T(w)=e^{T(w)},\,\,\mbox{ in $H^+$}. \end{equation}
Since the function
$$F(w)={w-1\over w+1}$$
is a conformal mapping of $H^+$ onto the unit disk ${\Bbb D},$
we see that the function
$$T(w)=\log{8|F'(w)|^2\over(1-|F(w)|^2)^2}=-2\log {\rm Re}\,w+\log 2$$
gives us a blow-up solution to the equation (\ref{eq4ab}) in $H^+.$
Now, by formula (\ref{factor5}), we have that the first required
solution has the form
$$u(z)=T(\omega)=-2\log{\rm Re}\,\omega(z)+\log 2=-2\log{\rm Re}\,x+\log 2.$$
\par
In order to get the second solution, we have to note that the function
$$F(w)=e^{-\lambda w},\,\,\lambda>0,$$
is also  conformal mapping of $H^+$ onto the punctured unit disk
${\Bbb D}\setminus\{0\}.$ Repeating the above arguments and taking
into account that ${\rm Re}\,\omega(z)=x,$ we arrive to solution
(\ref{eq1.32})
$$u(z)=T(\omega)=\log{8\lambda^2|e^{-\lambda\omega(z)}|^2\over (1-|e^{-\lambda\omega(z)}|^2)^2}=\log8\lambda^2-2 \lambda x-2\log(1-e^{-2\lambda x}).$$
\end{proof}

\par
  As an example, one can take the matrix function $A$ with the following constant entries $a_{11}=1,$ $a_{12}=-2,$ $a_{22}=5.$
Then $\mu(z)=(1+i)/2$ and $\omega(z)=x+i(y+2x).$ It is easy to
verify by straightforward computation that the functions
(\ref{eq1.31}), (\ref{eq1.32}) solve  the equation (\ref{eq1.30a}). In the
general case the admissible matrices in Corollary 5.3 have the
entries: $a_{11}=1,$ $a_{12}=\pm 2\nu(x)/\sqrt{1-\nu^2(x)},$
$a_{22}=(1+3\nu^2(x))/(1-\nu^2(x)).$ The particular case considered above corresponds  to
the function $ \nu(x)=1/\sqrt 2.$

\bigskip

Similar results can be given for the upper half-plane on the base of Corollary 5.4.

\par\medskip

\section{Free boundary}

The very important applied effect of the "dead zone" for solutions
of some partial differential equations, see e.g. \cite{Diaz}, the
Introduction and \S 1, is that a solution of the corresponding
differential equation vanishes on some nonempty open proper subset
of the domain of its definition. For example, it is well-known that
a solution of the semi-linear equation
$$\triangle u=u^q$$
may have the "dead zone" only when $0<q<1,$ see e.g. \cite{Diaz}, p.
15.
\par
We confine  ourselves to only one result in this respect, which is a
simple consequence of Proposition 5.3 and Corollary 4.2.

\bigskip

{\bf Theorem 6.1.} {\it Let ${\Bbb C}$ be  the complex plane
 and
 let the matrix function $A(z),$  $z=x+iy,$ be generated by  the Beltrami coefficient
\begin{equation}\label{eq1.60}
\mu(x)=\nu^2(x)\pm i\,\nu(x)\sqrt{1-\nu^2(x)},
\end{equation}
that is
\begin{equation}\label{matrix10}
A(z)=\left(\begin{array}{ccc} 1  & \mp {2\nu(x)\over\sqrt{1-\nu^2(x)}} \\
                             \mp {2\nu(x)\over\sqrt{1-\nu^2(x)} }         & {1+3\nu^2(x)\over 1-\nu^2(x)}
                              \end{array}\right),
\end{equation} where $\nu(x),$ $x\in{\Bbb R},$ stands for an arbitrary
measurable real-valued function, such that $|\nu(x)|\leq k<1.$ Then
the  semi-linear equation \begin{equation}\label{eq1.61} {\rm div\,}[A(z)\nabla
u]=u^q,\,\,\,0<q<1,\,\,z\in{\Bbb C}, \end{equation} has in the complex plane
the following solution with the "dead zone": \begin{equation}\label{matrix11}
u(x,y)=\left\{\begin{array}{ccc} &\gamma \left(y\pm\int\limits_0^x{2\nu(t)\over \sqrt{1-\nu^2(t)}}dt\right)^{2\over 1-q},
\,\,\,&if\,\, y>\varphi_{\pm}(x),\,x\in{\Bbb R},\\
                                 & 0                                                                                            &if\,\, x\leq \varphi_{\pm}(x).
                                 \end{array}\right.
                                 \end{equation}
 Here
$$\gamma=\left((1-q)^2\over 2(1+q)\right)^{1\over 1-q},$$
and
$$y=\varphi_{\pm}(x)=\pm\int\limits_0^x{2\nu(t)dt\over \sqrt{1-\nu^2(t)}},\,\,\-\infty<x<+\infty,$$
stands for the corresponding free boundary parametrization. }

\bigskip

\begin{proof}  Let $A(z)$ be a matrix function generated by the complex
Beltrami coefficient $\mu,$ satisfying (\ref{eq1.60}). Since $\mu$
depends on $x$ only, then formulae (\ref{eq1.11}), (\ref{eq1.12})
represent a unique quasiconformal mapping \begin{equation}\label{f1}
\omega(z)=x+i\left(y\pm  \int\limits_0^x{2\nu(t)\over
\sqrt{1-\nu^2(t)}}dt\right) \end{equation}
 of the complex plane ${\Bbb C}$ onto itself with complex dilatation $\mu$ normalize by
$\omega(0)=0,$ $\omega(i)=i,$ $\omega(\infty)=\infty.$
\par
 Since  $J_\omega(z)\equiv 1,$ one can apply Corollary 4.2 to represent
 solutions to the equation (\ref{eq1.61}) in the form
 $$u(z)=T(\omega(z)),$$
 where $T(w)$ satisfies
 the equation
 $$\triangle T(w)=T^q(w),\,\,w=\xi+i\eta.$$
 We see that the function
  $$T(w)=\gamma \eta^{2\over 1-q},\,\,\mbox{if\,\, $\eta>0$}$$
 and $T(w)=0,$ if $\eta\leq 0,$ satisfy the above  equation.
Thus we arrive at the required solution to the equation (\ref{eq1.61}):
$$u(x,y)=\gamma \left(y\pm\int\limits_0^x{2\nu(t)\over \sqrt{1-\nu^2(t)}}dt\right)^{2\over 1-q}, \,\,\,\mbox{if\,\, $y>\varphi_{\pm}(x),\,x\in{\Bbb R},$}
$$
and $u(z)=0,$ if $x\leq \varphi_{\pm}(x)$ where
$$\gamma=\left((1-q)^2\over 2(1+q)\right)^{1\over 1-q},\,\,\,\,\varphi_{\pm}(x)=\pm\int\limits_0^x{2\nu(t)dt\over \sqrt{1-\nu^2(t)}}.$$
\end{proof}

\par
 Let us verify this solution by straightforward computation. Since
$$u(x,y)=\gamma\left(y-\int\limits_0^x a_{12}(t)dt\right)^{2\over 1-q},$$
 we see that
$$u_x=-a_{12}\gamma{2\over 1-q}(y-\varphi_{+}(x))^{1+q\over 1-q},\,\, \mbox{and}\,\, u_y=\gamma{2\over 1-q}(y-\varphi_{+}(x))^{1+q\over 1-q}.$$
It yields that
$$a_{11}u_x+a_{12}u_y=-a_{12}\gamma{2\over 1-q}(y-\varphi_{+}(x))^{1+q\over 1-q}+a_{12}\gamma{2\over 1-q}(y-\varphi_{+}(x))^{1+q\over 1-q}=0.$$
On the other hand, since $a_{22}-a_{12}^2=1,$ we get
$$a_{12}u_x+a_{22}u_y=\gamma{2\over 1-q}(y-\varphi_{+}(x))^{1+q\over 1-q}.$$
Thus,
$${\rm div}\,(A(z)\nabla u)=(a_{12}u_x+a_{22}u_y)_y'=\gamma^q(y-\varphi_{+}(x))^{2q\over 1-q}=u^q,$$
and we complete the verification.

\section{Heat equation}

Factorization theorems, similar to Theorem 4.1 and its corollaries, also hold for parabolic and hyperbolic linear and semi-linear equations in the plane,
which
contain  the corresponding divergence operator in the linear part. We give, without proof, one  of such results.

\bigskip

{\bf Theorem 7.1.} {\it Let $\Omega$ be a domain in ${\Bbb C}$,
$A\in M^{2\times 2}(\Omega)$ and  let $f:{\Bbb R}\to{\Bbb R}$ be a
continuous function. Then  every weak solution of the semi-linear
equation \begin{equation}\label{qkk12a} 
u_t\, -\, a^2\,{\rm div\,} [A(z)\, \nabla u]\, =\, f(u) \ \ \ \mbox{in$\ \Omega\, ,\ t>0$} \, , 
\end{equation} can be
represented as the composition \begin{equation}\label{factork} u(z,t)\, =\, T(\,
\omega(z),\, t\, )\, , \end{equation} where $\omega:\Omega\to G$ is a
quasiconformal mapping agreed with $A$, and $T$ is a weak solution
of the equation \begin{equation}\label{eq4ak} T_t\, -\, a^2\,J^{-1}
\triangle\,T\, =\, f(T)\ \ \ \ \mbox{in\ $G$\ ,\ $t>0$}\ . \end{equation} Here
$J$ stands for the Jacobian of the mapping $\omega.$ }

\bigskip

As an example for the illustration of Theorem 6.1. we give a statement concerning the heat
conductivity equation in  divergent form.

\bigskip

{\bf Proposition 7.1.} {\it  Let ${\Bbb D}$ be  the unit disk and
 let the matrix function $A(z)$ be generated by  the Beltrami coefficient
\begin{equation}\label{fe3} \mu(z)=k(|z|){z\over\bar z} \end{equation} in accordance with
formula (\ref{matrix}) where \begin{equation}\label{fe4} k(t)=\nu^2(t)\pm
i\,\nu(t)\sqrt{1-\nu^2(t)} \end{equation} and $\nu(t),$ $0\leq t<1,$ stands for
an arbitrary measurable real-valued function, such that
$|\nu(t)|\leq q<1.$ Then the fundamental solution to the
inhomogeneous anisotropic  heat equation
\begin{equation}\label{fe5}
 u_t-a^2\,{\rm div}\,(A(z)\nabla u)=0,\,\,z\in {\Bbb D},\,\,t>0,
\end{equation}
 has the representation
\begin{equation}\label{fe6}
 u(z,t)={1\over 4\pi a^2t}\exp\left(-{|z|^2\over 4a^2t}\right),\,\,z\in {\Bbb D},\,\, t>0.
\end{equation}
 }

\bigskip

 \begin{proof} We are looking for a solution in the form $u(z,t)=T(\omega(z),t),$ where $\omega(z)$ stands for a quasiconformal automorphism
 of the unit disk $\mathbb D$  agreed with the matrix function $A(z).$ Since the mapping $w=\omega(z)$ is volume preserving, we see that its Jacobian is identically equal to the unit.
   Then, by Theorem 7.1, $T(w,t)$ will satisfy the equation
\begin{equation}\label{fe7}
 T_t-a^2\triangle T=0,\,\,w\in {\Bbb D},\,\, t>0.
\end{equation}
 It is a classical result that the
 fundamental solution to the inhomogeneous heat equation (\ref{fe7})
 has the representation
\begin{equation}\label{fe8}
 T(w,t)={1\over(2a\sqrt{\pi t})^2}
\exp\left(-{|w|^2\over 4a^2t}\right),\,\,w\in {\Bbb D},\,\, t>0.
\end{equation}
 Since in our case $|\omega(z)|=|z|,$ we arrive at the required conclusion.
\end{proof}

\bigskip

 {\bf Remark 7.1.} In a similar way we can study the well-known semi-linear combustion equation  in an anisotropic medium
\begin{equation}\label{eq7.1} u_t\, -\, a^2\,{\rm div\,}[A(z)\, \nabla u]\, =\,
e^u \end{equation} by means of reducing its weak solution to the composition of
a weak solution in the isotropic case, see e.g. \cite{Barenblat} and
\cite{Pohozhaev} and the references therein, \begin{equation}\label{eq7.2} T_t\,
-\, a^2\,J^{-1}\triangle\,T\, =\, e^T \end{equation} with a quasiconformal
mapping $\omega $ generated by $A$. As above, $J$ stands for the
Jacobian of the mapping $\omega .$

 \section{Boundary value problems}

 The Factorization Theorem 4.1 allows us, in particular, to reduce the study of the boundary value problems, say the Dirichlet problem
\begin{equation}\label{matrix21}
\left\{\begin{array}{ccc} &{\rm div}\,(A(z)\nabla u)=f(u) & \textrm{if}\,\, z\in\Omega,\\
                          & u=\varphi,                             & \textrm{if}\,\, z\in\partial\Omega
                                 \end{array}\right.
                                 \end{equation}
in a Jordan domain $\Omega$ with $A\in M^{2\times 2}(\Omega)$ and
continuous boundary function $\varphi$ to the study of the following
Dirichlet problem in the unit disk \begin{equation}\label{matrix22}
\left\{\begin{array}{ccc} & \triangle T=J(w)f(T) &if\,\, |w|<1,\\
                          & T=\psi ,                             &if\,\, |w|=1
                                 \end{array}\right.
                                 \end{equation}
with the corresponding boundary function
$\psi=\varphi\circ\omega^{-1}$ where $\omega $ is a quasiconformal
mapping of $\Omega$ onto the unit disk $\mathbb D$ agreed with $A$,
and $J$ is the Jacobian of $\omega^{-1}$. First of all, note that
$\omega$ exists by the Measurable Riemann mapping theorem, see e.g.
Theorem V.B.3 in \cite{Ahlfors:book} and Theorem V.1.3 in
\cite{LV:book}.

The existence and continuity of the boundary function $\psi$ in the
case of an arbitrary Jordan domain $\Omega$ is a fundamental result
of the theory of the boundary behavior of conformal and
quasiconformal mappings. Namely,
$$\omega^{-1}=h\circ g$$
where $g$ stands for a quasiconformal automorphism of the unit disk
${\mathbb{D}}$ and $h$ is a conformal mapping of $\mathbb{D}$ onto
$\Omega.$ It is known that $g$ can be extended to a homeomorphism of
$\overline{\Bbb D}$ onto itself, see e.g. Theorem I.8.2 in
\cite{LV:book}. Moreover, by the well-known
Caratheodory-Osgood-Taylor theorem on the boundary correspondence
under conformal mappings, see \cite{C} and \cite{OT}, the mapping
$h$ is extended to a homeomorphism of $\overline{\Bbb D}$ onto
$\overline{\Omega}$. Thus, the function $\psi$ is well defined and
really continuous on the unit circle.

\noindent {\it Institute of Applied Mathematics and Mechanics of
National Academy\\ of Sciences of Ukraine, 84100, Ukraine,
Slavyansk, Dobrovolskogo str. 1},

\noindent Email: vgutlyanskii@gmail.com, star-o@ukr.net,
vlryazanov1@rambler.ru

\end{document}